\documentclass[11pt,a4paper,leqno]{article}
\usepackage{a4wide}
\usepackage{latexsym}
\usepackage{amsmath}
\usepackage{amsthm}
\usepackage{amssymb}
\theoremstyle{plain}
\newtheorem{theo}{Theorem}[section]
\newtheorem{lemma}[theo]{Lemma}
\newtheorem{prop}[theo]{Proposition}

\theoremstyle{definition}
\newtheorem{defin}[theo]{Definition}
\newtheorem{rem}[theo]{Remark}

\newenvironment{proof1}{\medskip\par\noindent{\bf Proof}}{\hfill $\Box$
\medskip\par}

\def\C{\mathbb{C}}
\def\N{\mathbb{N}}
\def\R{\mathbb{R}}

\def\a{\alpha}
\def\b{\beta}
\def\ga{\gamma}
\def\balpha{\boldsymbol{\alpha}}
\def\bbeta{\boldsymbol{\beta}}
\def\bdelta{\boldsymbol{\delta}}
\def\blambda{\boldsymbol{\lambda}}
\def\bro{\boldsymbol{\rho}}
\def\btheta{\boldsymbol{\theta}}

\def\ba{\boldsymbol{a}}
\def\bd{\boldsymbol{d}}
\def\be{\boldsymbol{e}}
\def\bfe{\boldsymbol{f}}
\def\bJ{\boldsymbol{J}}
\def\bL{\boldsymbol{L}}
\def\bM{\boldsymbol{M}}
\def\bm{\boldsymbol{m}}
\def\bN{\boldsymbol{N}}
\def\bw{\boldsymbol{w}}
\def\bz{\boldsymbol{z}}

\begin{document}

\title{Continuous right inverses for the asymptotic Borel map in ultraholomorphic classes via a Laplace-type transform}
\author{Alberto Lastra, St\'ephane Malek and Javier Sanz}
\date{\today}

\maketitle

{ \small
\begin{center}
{\bf Abstract}
\end{center}
 A new construction of linear continuous right inverses for the asymptotic Borel map is provided in the framework of general Carleman ultraholomorphic classes in narrow sectors. Such operators were already obtained by V. Thilliez by means of Whitney extension results for non quasianalytic ultradifferentiable classes, due to J. Chaumat and A. M. Chollet, but our approach is completely different, resting on the introduction of a suitable truncated Laplace-type transform. This technique is better suited for a generalization of these results to the several variables setting. Moreover, it closely resembles the classical procedure in the case of Gevrey classes, so indicating the way for the introduction of a concept of summability which generalizes $k-$summability theory as developed by J. P. Ramis.

}

\bigskip \bigskip

\section{Introduction}

\indent\indent For sectors $S$ of suitably small opening and vertex at 0, the Borel-Ritt-Gevrey theorem proved by J. P. Ramis (see~\cite{Ramis1,Ramis2,MartinetRamis}, \cite[Thm.\ 2.2.1]{Balser}) guarantees the existence of holomorphic functions on $S$ having an arbitrarily prescribed Gevrey asymptotic expansion of order $\alpha$ at~0.
This amounts to the surjectivity of the asymptotic Borel map, sending a function to its series of asymptotic expansion,
when considered between the corresponding spaces of Gevrey functions, respectively Gevrey series. The proof is constructive,
and basically consists in applying a truncated Laplace transform to the formal Borel transform of the initially given Gevrey series.

For functions $f$ holomorphic on a polysector $S\subset\C^n$ with vertex at $\bf 0$, H.~Majima~\cite{Majima1,Majima2} put forward
the concept of strong asymptotic developability, which has been shown~\cite{hernandez,galindosanz} to amount to the boundedness of the
derivatives of $f$ on bounded proper subpolysectors of $S$, just as in the one-variable situation.
The asymptotic behaviour of $f$ is determined by the family $\mathrm{TA}(f)$ (see Section~\ref{seccSeverVaria}),
consisting of functions obtained as limits of the derivatives of $f$
when some of its variables tend to 0 (in the same way as the coefficients of the series of asymptotic expansion in the one-variable case).

In 1989 Y. Haraoka~\cite{Haraoka} considered the space of holomorphic functions $f$ in a polysector $S$ that admit Gevrey strong asymptotic expansion of order
$\balpha=(\a_1,\ldots,\a_n)\in[1,\infty)^n$ (one order per variable), and got two partial Borel-Ritt-Gevrey type results in this context again by applying a
(multidimensional) truncated Laplace transform.

Subsequently, in the one-variable setting V.~Thilliez~\cite[Theorem 1.3]{Thilliez1} obtained a linear continuous version of this result by constructing
extension operators (linear continuous right inverses for the Borel map) from Banach spaces of Gevrey series into Banach
spaces of functions whose derivatives admit Gevrey-like bounds uniformly on all of $S$ (so that they admit Gevrey asymptotic expansion at 0). His proof rests
on Whitney type extension results for ultradifferentiable classes by J.~Chaumat and A.~M.~Chollet~\cite{ChaumatChollet}.

The third author of the present work re-proved in~\cite{javier} Thilliez's result
in an elementary way by a careful study of Ramis' argument. The solution so obtained, in integral form,  is valid for vector, Banach space-valued functions, and it is also amenable to the determination of its
behaviour in case this Banach space consists precisely of Gevrey functions. Since these Banach spaces verify an exponential-law isomorphism, one may apply a recurrent argument on the number of variables to obtain extension operators in several variables which generalize Thilliez's result and provide linear continuous versions of the first interpolation result proven by Haraoka~\cite[Theorem 1.(1)]{Haraoka} and a right inverse for the map $f\mapsto\mathrm{TA}(f)$.

The next step in these developments was again taken by V. Thilliez in~\cite{thilliez}, where he broadens the scope of the preceding one-dimensional results on considering general ultraholomorphic classes in sectors.
Specifically, given $A>0$, a sequence of positive real numbers $\bM=(M_{p})_{p\in\N_{0}}$ and a sector $S$ with vertex at $0$ in the Riemann surface of the logarithm, $\mathcal{R}$,
$\mathcal{A}_{\bM,A}(S)$ consists of the complex holomorphic functions $f$ defined in $S$ such that
$$\Vert f\Vert_{\bM,A,S}:=\sup_{p\in\N_{0},z\in S}\frac{|D^{p}f(z)|}{A^{p}p!M_{p}}<\infty.$$
The ultraholomorphic Carleman class $\mathcal{A}_{\bM}(S)$ is defined as $\cup_{A>0}\mathcal{A}_{\bM,A}(S)$.\par
Accordingly, $\Lambda_{\bM,A}(\N_{0})$ is the set of the sequences of complex numbers
$\blambda=(\lambda_{p})_{p\in\N_{0}}$ such that
$$|\blambda|_{\bM,A}:=\sup_{p\in\N_{0}}\frac{|\lambda_{p}|}{A^{p}p!M_{p}}<\infty,$$
and $\Lambda_{\bM}(\N_0):=\cup_{A>0}\Lambda_{\bM,A}(\N_0)$.
$(\mathcal{A}_{\bM,A}(S),\Vert \cdot\Vert_{\bM,A,S})$ and $(\Lambda_{\bM,A}(\N_{0}),|\cdot|_{\bM,A})$ are Banach spaces,
and, as the derivatives of the elements in $\mathcal{A}_{\bM,A}(S)$ are Lipschitzian, we may define the (linear and continuous) asymptotic Borel map $\mathcal{B}:\mathcal{A}_{\bM,A}(S)\to\Lambda_{\bM,A}(\N_{0})$ given by
$$
\mathcal{B}(f):=
\big(f^{(p)}(0)\big)_{p\in\N_{0}}\in\C^{\N_{0}},\qquad f^{(p)}(0):=\lim_{z\to0}f^{(p)}(z).
$$
Gevrey classes of order $\a>1$ in a sector $S$ correspond to the sequence $\bM_{\a}=(p!^{\a-1})_{p\in\N_0}$.
For strongly regular sequences $\bM$ (see Subsection~\ref{sucfuereg}), among which we find the sequences $\bM_{\a}$, the construction of Thilliez's operators in the next theorem is based on a double application of suitable Whitney's extension results for Whitney ultradifferentiable jets on compact sets with Lipschitz boundary, given by
J.~Chaumat and A.~M.~Chollet in~\cite{chaucho}, and on a solution of a $\overline{\partial}$-problem.

\begin{theo}[\cite{thilliez},\ Thm.\ 3.2.1]
\label{teorthilliez}
Let $\bM=(M_{p})_{p\in\N_{0}}$ be a strongly regular sequence with associated growth index $\ga(\bM)$. Let us consider $\ga\in\R$ with $0<\ga<\ga(\bM)$,
and let $S_{\gamma}$ be a sector with opening $\ga\pi$. Then there exists $d\ge 1$, that only depends on $\bM$ and $\ga$, so that for every $A>0$ there exists a linear continuous operator
$$T_{\bM,A,\ga}:\Lambda_{\bM,A}(\N_{0})\longrightarrow \mathcal{A}_{\bM,dA}(S_{\ga})$$
such that $\mathcal{B}\circ T_{\bM,A,\ga}\blambda=\blambda$ for every $\blambda\in\Lambda_{\bM,A}(\N_{0})$.
\end{theo}

For Gevrey classes, $\ga(\bM_{\a})=\a-1$, so that the condition in the theorem tells that the opening of the sector should be less than $(\a-1)\pi$ for the extension to exist,
what agrees with the classical Borel-Ritt-Gevrey statement.

This result has been extended to functions of several variables by the first and third authors~\cite{lastrasanz} by applying a recursive technique similar to that in~\cite{javier}, but resting on this new construction of Thilliez, what makes it difficult to determine the behaviour of the derivatives of the solution of the one dimensional problem when it takes its values in a Banach space of the type $\mathcal{A}_{\bM,A}(S)$. As indicated above, this information is crucial in the process providing a right inverse for the map $\mathrm{TA}$ in this context.

With these preliminaries, the main aim in the present work is to obtain a new proof of Theorem~\ref{teorthilliez} which no longer depends on Whitney-type extension results, but rather makes use of a suitable truncated integral, Laplace-like operator, in the same vein as Ramis' original proof. The kernel in this integral operator will be given in terms of a flat function obtained by V. Thilliez~\cite[Thm.\ 2.3.1]{thilliez}, playing a similar role as that played by the exponential $\exp(-1/z^{1/\a})$ in the Gevrey case of order $\a$. Indeed, in the authors' opinion the absence of an elementary function governing null asymptotics in this general case was the reason for the use, up to this moment, of results belonging to the ultradifferentiable setting when solving interpolation problems in non-Gevrey ultraholomorphic classes. As stated before, this new approach is better suited for the generalization to the several variables setting, and moreover, it provides some insight when searching for a summability tool in general ultraholomorphic classes which resembles $k-$summability, specifically designed for the Gevrey case and which has proved itself extremely useful in the reconstruction of analytic solutions of linear and nonlinear (systems of) meromorphic ordinary differential equations at irregular singular points, departing from their formal power series solutions (see~\cite{balser} and the references therein). We include in the last section some hints in this direction, where we will make use of quasi-analyticity properties in these classes which have been characterized (see~\cite{lastrasanz1}) in terms of Watson's type lemmas. It should also be indicated that the construction of the formal and analytic transforms incorporated into this new technique is inspired by the study of general summability methods, equivalent in a sense to $k-$summability, developed by W.~Balser in~\cite[Section\ 5.5]{balser} and which have already found its application to the analysis of formal power series solutions of different classes of partial differential equations and so-called moment-partial differential equations (see the works of W. Balser and Y. Yoshino~\cite{BalserYoshino}, the second author~\cite{Malek1,Malek2} and S. Michalik~\cite{Michalik1,Michalik2,Michalik3,Michalik4}, among others).
Also, some results on summability for non-Gevrey classes, associated to strongly regular sequences, have been provided for difference equations
by G. K. Immink in~\cite{Immink}, whereas V. Thilliez has obtained some results on solutions within these general classes for algebraic equations in~\cite{Thilliez2}.
We hope our summability theory is able to shed some light on some of these problems or on similar ones.

\section{Preliminaries}
\subsection{Notation}\label{notation}
We set $\N:=\{1,2,...\}$, $\N_{0}:=\N\cup\{0\}$.
$\mathcal{R}$ stands for the Riemann surface of the logarithm, and
$\C[[z]]$ is the space of formal power series in $z$ with complex coefficients.\par\noindent
For $\gamma>0$, we consider unbounded sectors
$$S_{\gamma}:=\{z\in\mathcal{R}:|\hbox{arg}(z)|<\frac{\gamma\,\pi}{2}\}$$
or, in general, bounded or unbounded sectors
$$S(d,\alpha,r):=\{z\in\mathcal{R}:|\hbox{arg}(z)-d|<\frac{\alpha\,\pi}{2},\ |z|<r\},\quad
S(d,\alpha):=\{z\in\mathcal{R}:|\hbox{arg}(z)-d|<\frac{\alpha\,\pi}{2}\}$$
with bisecting direction $d\in\R$, opening $\alpha\,\pi$ and (in the first case) radius $r\in(0,\infty)$.\par\noindent
A sectorial region $G(d,\a)$ will be a domain in $\mathcal{R}$ such that $G(d,\a)\subset S(d,\a)$, and
for every $\beta\in(0,\a)$ there exists $\rho=\rho(\beta)>0$ with $S(d,\beta,\rho)\subset G(d,\a)$.\par\noindent
$D(z_0,r)$ stands for the disk centered at $z_0$ with radius $r>0$.

\noindent For $n\in\N$, we put $\mathcal{N}=\{1,2,\ldots,n\}$. If $J$ is a nonempty subset of $\mathcal{N}$,  $\#J$ denotes its cardinal number.\par\noindent
A polysector is a product of sectors, $S=\prod_{j=1}^n S_j\subset\mathcal{R}^n$.
The polysector $\prod_{j=1}^n S(d_j,\theta_j,\rho_j)$ (with $\rho_j$ possibly equal to $\infty$)
will be denoted by $S=S(\bd, \btheta, \bro)$, with the obvious meaning for $\bd$, $\btheta$ and $\bro$. In case $\rho_j=+\infty$ for
$j\in\mathcal{N}$, we write $S=S(\bd,\btheta)$.\par\noindent
We say a polysector $T=\prod_{j=1}^n T(d'_j,\theta'_j,\rho'_j)$
is a {\it bounded proper subpolysector}
of~$S=S(\bd, \btheta, \bro)$, and we write $T\ll S$, if
for $j\in\mathcal{N}$ we have $\rho'_j<\rho_j$ (so that $\rho'_j<+\infty$) and
\begin{equation}\label{subsectorpropio}
[d'_j-\theta'_j/2,d'_j+\theta'_j/2]\subset (d_j-\theta_j/2,d_j+\theta_j/2).
\end{equation}
Finally, we say $T=\prod_{j=1}^n T(d'_j,\theta'_j)$
is an {\it unbounded proper subpolysector}
of~$S=S(\bd, \btheta)$, and we write $T\prec S$, if
for $j\in\mathcal{N}$ we have $(\ref{subsectorpropio})$.
Given $\bz\in\mathcal{R}^n$, we write $\bz_J$ for the restriction of $\bz$ to $J$, regarding
$\bz$ as an element of $\mathcal{R}^{\mathcal{N}}$.\par\noindent
Let $J$ and $L$ be nonempty disjoint subsets of $\mathcal{N}$. For
$\bz_J\in\mathcal{R}^J$ and $\bz_L\in\mathcal{R}^L$, $(\bz_J,
\bz_L)$ represents the element of $\mathcal{R}^{J\cup L}$ satisfying
$
(\bz_J, \bz_L)_J=\bz_J
$,
$
(\bz_J, \bz_L)_L=\bz_L
$;
we also write
$J^{\prime}=\mathcal{N}\setminus J$, and for $j\in \mathcal{N}$ we use $j^{\prime}$ instead of $\{j\}^{\prime}$.
In particular, we shall use these conventions for
multi-indices.\par\noindent
For $\btheta=(\theta_1,\ldots,\theta_n)\in(0,\infty)^n$, we write $S_{\btheta}=\prod_{j=1}^nS_{\theta_j}$ and $S_{\btheta_J}=\prod_{j\in J}S_{\theta_j}\subset\mathcal{R}^J$.\par\noindent
If $\boldsymbol{z}=(z_{1},z_{2},\ldots,z_{n})\in\mathcal{R}^{n}$, $\balpha=(\a_{1},\a_{2},\ldots,\a_{n})$, $\bbeta=(\b_{1},\b_{2},\ldots,\b_{n})\in\N_0^{n}$, we define:
$$\begin{array}{lll}
|\balpha|=\a_{1}+\a_{2}+\ldots+\a_{n}, & \balpha!=\a_{1}!\a_{2}!\cdots\a_{n}!, \\
D^{\balpha}=\frac{\partial^{\balpha}}{\partial\boldsymbol{z}^{\balpha}}=
\frac{\partial^{|\balpha|}}{\partial z_{1}^{\a_{1}}\partial z_{2}^{\a_{2}}\ldots\partial z_{n}^{\a_{n}}}, & \be_j=(0,\ldots,\stackrel{j)}{1},\ldots,0).
\end{array}$$
For $\bJ\in\N_0^n$, we will frequently write $j=|\bJ|$.

\subsection{Asymptotic expansions}

Given $A>0$, a sequence of positive real numbers $\bM=(M_{p})_{p\in\N_{0}}$ and a sector $S$, for every
$f$ in the class $\mathcal{A}_{\bM,A}(S)$ one may put
$$f^{(p)}(0):=\lim_{z\in S,z\to0 }f^{(p)}(z)\in\C$$
for every $p\in\N_{0}$.
Then, $f$ admits the formal power series $\sum_{p\in\N_0}\frac{1}{p!}f^{(p)}(0)z^p$ as its uniform asymptotic expansion at 0, in the following sense.
\begin{defin}
Let $\bM=(M_{p})_{p\in\N_0}$ be a sequence of positive real numbers and let $f$ be a holomorphic function in a sector $S$ with vertex at the origin. We say $f$ admits the formal power series $\hat{f}=\sum_{p=0}^{\infty}a_{p}z^{p}\in\C[[z]]$ as its \textit{uniform $\bM-$asymptotic expansion} in $S$ of type $A>0$ (when the variable tends to 0) if there exists $C>0$ such that for every $N\in\N$, one has
\begin{equation}\left|f(z)-\sum_{p=0}^{N-1}a_pz^p \right|\le CA^NM_{N},\qquad z\in S.\label{desarasintunifo}
\end{equation}
We will write $f\sim_{\bM}\sum_{p=0}^{\infty}a_pz^p$ (uniformly in $S$ and with type $A$).
\end{defin}

\begin{rem}\label{remaCarlclassasympexpan}
Conversely, and as a consequence of Cauchy's integral formula for the derivatives, one can prove that whenever $T$ is a proper subsector of $S$, there exists a constant $c=c(T,S)>0$ such that the restriction to $T$, $f_T$, of functions $f$ defined on $S$ and admitting uniform $\bM-$asymptotic expansion in $S$ of type $A>0$, belongs to $\mathcal{A}_{\bM,cA}(T)$, and moreover, if one has (\ref{desarasintunifo}) then $\Vert f_T\Vert_{\bM,cA,T}\le C$.
\end{rem}

\begin{rem}\label{remaasympexpansectoregio}
For sectorial regions $G$, $f\sim_{\bM}\sum_{p=0}^{\infty}a_pz^p$ in $G$ means that $f\sim_{\bM}\sum_{p=0}^{\infty}a_pz^p$ uniformly in every sector $S$ such that $\overline{S}\setminus\{0\}\subset G$.
\end{rem}

\subsection{Strongly regular sequences}\label{strregseq}

The information in this subsection is taken from the work of V. Thilliez~\cite{thilliez}, which we refer to for further details and proofs.
In what follows, $\bM=(M_p)_{p\in\N_0}$ will always stand for a sequence of
positive real numbers, and we will always assume that $M_0=1$.
\begin{defin}\label{sucfuereg}
We say $\bM$ is \textit{strongly regular} if the following hold:\par
($\a_0$) $\bM$ is \textit{logarithmically convex}: $M_{p}^{2}\le M_{p-1}M_{p+1}$ for
every $p\in\N$.\par
($\mu$) $\bM$ is of \textit{moderate growth}: there exists $A>0$ such that
$$M_{p+\ell}\le A^{p+\ell}M_{p}M_{\ell},\qquad p,\ell\in\N_0.$$
\par($\gamma_1$) $\bM$ satisfies the \textit{strong non-quasianalyticity condition}: there exists $B>0$ such that
$$
\sum_{\ell\ge p}\frac{M_{\ell}}{(\ell+1)M_{\ell+1}}\le B\frac{M_{p}}{M_{p+1}},\qquad p\in\N_0.
$$
\end{defin}
For a strongly regular sequence $\bM=(M_{p})_{p\in\N_0}$, it is direct to check from properties $(\alpha_0)$ and $(\gamma_{1})$ that $\bm=(m_{p}:=M_{p+1}/M_{p})_{p\in\N_0}$ is an increasing sequence to infinity, so that the map $h_{\bM}:[0,\infty)\to\R$, defined by
$$h_{\bM}(t):=\inf_{p\in\N_{0}}M_{p}t^p,\qquad h_{\bM}(0)=0$$
turns out to be a non-decreasing continuous map in $[0,\infty)$, and its range is the set $[0,1]$. In fact
$$
h_{\bM}(t)= \left \{ \begin{matrix}  t^{p}M_{p} & \mbox{if }t\in\left[\frac{1}{m_{p}}\right. ,\left. \frac{1}{m_{p-1}}\right),\ p=1,2,\ldots,\\
1 & \mbox{if } t\ge 1/m_{0}. \end{matrix}\right.
$$
Some properties of strongly regular sequences needed in the present work are the following.
\begin{lemma}
Let $\bM=(M_{p})_{p\in\N_0}$ be a strongly regular sequence and $A>0$ the constant appearing in $(\mu)$. Then,
\begin{align}
\label{propgrowth}
&M_{p+\ell}\ge M_pM_{\ell},\qquad \hbox{for every }p,\ell\in\N_0,\\
\label{e107}
&m_{p}\le A^{2}M_{p}^{1/p},\qquad \hbox{for every }p\in\N_0,\\
\label{e115}
&M_{p}^{1/p}\le m_{p},\qquad \hbox{for every }p\in\N_0.
\end{align}
Let $s$ be a real number with $s\ge1$. There exists $\rho(s)\ge1$ (only depending on $s$ and $\bM$) such that
\begin{equation}\label{e120}
h_{\bM}(t)\le(h_{\bM}(\rho(s)t))^{s}\qquad\hbox{for }t\ge0.
\end{equation}
\end{lemma}

\begin{defin}
Let $\bM=(M_{p})_{p\in\N_{0}}$ be a strongly regular sequence, $\ga>0$. We say $\bM$ satisfies property $\left(P_{\ga}\right)$  if there exist a sequence of real numbers $m'=(m'_{p})_{p\in\N_0}$ and a constant $a\ge1$ such that: (i) $a^{-1}m_{p}\le m'_{p}\le am_{p}$, $p\in\N$, and (ii) $\left((p+1)^{-\ga}m'_{p}\right)_{p\in\N_0}$ is increasing.

The \textit{growth index} of $\bM$ is
$$\ga(\bM):=\sup\{\ga\in\R:(P_{\ga})\hbox{ is fulfilled}\}.$$
\end{defin}

For any strongly regular sequence $\bM$ one has $\gamma(\bM)\in(0,\infty)$.
For the Gevrey sequence of order $\a>0$ given by $\bM_{\a}=(p!^{\a})_{p\in\N_{0}}$, we have $\ga(\bM_{\a})=\a$.

Finally we describe the properties of a function that will be crucial in the construction of a kernel for our Laplace-type operator.
\begin{prop}[\cite{thilliez},\ Thm.\ 2.3.1 and Lemma\ 2.3.2]\label{gm}
Suppose $\bM=(M_{p})_{p\in\N_0}$ is a strongly regular sequence and $\delta\in\R$ with $0<\delta<\gamma(\bM)$. There exists a holomorphic function $G_{\bM}$ defined in $S_{\delta}$ such that for every $w\in S_{\delta}$ one has:
\begin{itemize}
\item[(i)] $k_{1}h_{\bM}(k_2|w|)\le|G_{\bM}(w)|\le h_{\bM}(k_3|w|)$, where $k_1$, $k_2$ and $k_3$ are positive constants that only depend on $\bM$ and $\delta$.
\item[(ii)] For every $p\in\N_0$, $|G_{\bM}^{(p)}(w)|\le b_1^pp!M_ph_{\bM}(b_2|w|)$, $b_1$ and $b_2$ being positive constants that only depend on $\bM$ and $\delta$. In particular, we deduce that  $G_{\bM}\in\mathcal{A}_{\bM}(S_{\delta})$ and it is flat, i.e., $G_{\bM}\sim_{\bM}0$ uniformly in $S_\delta$.
\item[(iii)]For every $p\in\N_0$, $|(1/G_{\bM})^{(p)}(w)|\le b_3b_4^pp!M_p(h_{\bM}(b_5|w|))^{-1}$, where $b_3$, $b_4$ and $b_5$ are positive constants that only depend on $\bM$ and $\delta$.
\end{itemize}
\end{prop}

\begin{rem}\label{remaConstrG_M}
Let $0<\delta<\gamma(\bM)$. The function $G_{\bM}$ is defined as follows. Take $\delta_1$ and $s$ with $\delta<\delta_1<\gamma(\bM)$ and $s\delta_1<1<s\gamma(\bM)$. Then
\begin{equation}\label{gthilliez} G_{\bM}(z)=\exp\left(\frac{1}{\pi}\int_{-\infty}^{\infty}\log\left(h_{\bM^s}(|t|)\right)\frac{itz^s-1}{it-z^s}\frac{dt}{1+t^2}\right),\qquad z\in S_{\delta_{1}},
\end{equation}
with $\bM^s:=(M_p^{s})_{p\in\N_0}$, which turns out to be a strongly regular sequence too.
The restriction of $G_{\bM}$ to $S_{\delta}$ is the function in Proposition~\ref{gm}.
\end{rem}

\section{Moment sequence associated to $\bM$}\label{seccion2}

This section is devoted to the construction of a moment function $e_{\bM}$, associated to a strongly regular sequence $\bM$, which in turn will provide us with a sequence of moments $\mathfrak{m}=(m(p))_{p\in\N_{0}}$ equivalent, in the sense of the following definition, to $\bM$.
\begin{defin}[see \cite{pet},~\cite{chaucho}]
Two sequences $\bM=(M_{p})_{p\in\N_0}$ and $\bM'=(M'_{p})_{p\in\N_0}$ of positive real numbers are said to be \textit{equivalent} if there exist positive constants $L,H$ such that
$$L^pM_p\le M'_p\le H^pM_p,\qquad p\in\N_0.$$
\end{defin}

We note that, given a sector $S$ and a pair of equivalent sequences $\bM$ and $\bM'$, the spaces $\mathcal{A}_{\bM}(S)$ and $\mathcal{A}_{\bM'}(S)$ coincide.

Let $\bM=(M_{p})_{p\in\N_0}$ be a strongly regular sequence with growth index $\gamma(\bM)$. We take $0<\delta<\gamma(\bM)$ and define $e_{\bM}:S_{\delta}\to\C$ by
\begin{equation}\label{funcione}
e_{\bM}(z):=zG_{\bM}(1/z), \qquad z\in S_{\delta},
\end{equation}
where $G_{\bM}$ is defined in Subsection~\ref{strregseq}.

\begin{rem}
There is some freedom in the choice of $e_{\bM}$. Firstly, the factor $z$ may be changed into any $z^{\alpha}$ for some positive real number $\a$ (so that the assertion $(i)$ in the next lemma holds true), where the principal branch of the power is considered. Our choice tries to make the following computations simpler. Secondly, as indicated in Remark~\ref{remaConstrG_M}, there are some constants $\delta_1$ and $s$ to be fixed in the construction of $G_{\bM}$.
\end{rem}

\begin{lemma}\label{propiedadese}
The function $e_{\bM}$ satisfies the following assertions:
\begin{itemize}
\item[(i)] $e_{\bM}$ is well defined in $S_{\delta}$ and is such that $z^{-1}e_{\bM}(z)$ is integrable at the origin, it is to say, for any $t_0>0$ and $\tau\in\R$ with $|\tau|<\frac{\delta\pi}{2}$ the integral $\int_{0}^{t_0}t^{-1}|e_{\bM}(te^{i\tau})|dt$ is finite.
\item[(ii)] There exist $C,K>0$ (not depending on $\delta$) such that
\begin{equation}\label{e147}|e_{\bM}(z)|\le Ch_{\bM}\left(\frac{K}{|z|}\right),\qquad z\in S_{\delta}.
\end{equation}
\item[(iii)] For $x\in\R$, $x>0$, the values of $e_{\bM}(x)$ are positive real.
\end{itemize}
\end{lemma}
\begin{proof1}
Let $t_0>0$ and $\tau\in\R$ with $|\tau|<\frac{\delta\pi}{2}$. From Proposition~\ref{gm} there exists $k_3>0$ such that
$$\int_{0}^{t_0}\frac{|e_{\bM}(te^{i\tau})|}{t}dt\le \int_{0}^{t_0}h_{\bM}(k_3/t)dt.$$
We conclude the convergence of the last integral from the fact that $h_{\bM}(s)\equiv 1$ when $s\ge\frac{1}{m_{1}}$ and its continuity in $[0,\infty)$. The first part of the result is achieved.\par
For the second, we have
$$|e_{\bM}(z)|=|z||G_{\bM}(1/z)|\le|z|h_{\bM}(k_3/|z|),$$
for every $z\in S_{\delta}$, so $(ii)$ holds for $|z|<\tilde{M}$ for any fixed $\tilde{M}>0$. If $|z|\ge \tilde{M}$, we apply (\ref{e120}) for $s=2$ and the very definition of $h_{\bM}$ to get
$$
|e_{\bM}(z)|\le |z|\Big(h_{\bM}\big(\frac{\rho(2)k_3}{|z|}\big)\Big)^2\le |z|h_{\bM}\big(\frac{\rho(2)k_3}{|z|}\big)M_2\big(\frac{\rho(2)k_3}{|z|}\big)^2\le \frac{\rho(2)^2k_3^2M_2}{\tilde{M}}h_{\bM}\big(\frac{\rho(2)k_3}{|z|}\big).
$$
Finally, if $x>0$, then $e_{\bM}(x)=xG_{\bM}(1/x)$. From (\ref{gthilliez}) we have
$$G_{\bM}(1/x)=\exp\left(\frac{1}{\pi}\int_{-\infty}^{\infty}\log\left(h_{\bN}(|t|)\right)\frac{it-x^s}{itx^s-1}\frac{dt}{1+t^2}\right).$$
It is immediate to check that the imaginary part of the expression inside the previous integral is odd with respect to $t$, so that the corresponding integral is 0 and $G_{\bM}(1/x)$ is positive and real for $x>0$.
\end{proof1}

The role that Eulerian Gamma function played for Gevrey sequences will now be played by the following auxiliary function.

\begin{defin}
We define the \textit{moment function} associated to $\bM$ as
$$m(\lambda):=\int_{0}^{\infty}t^{\lambda-1}e_{\bM}(t)dt=\int_{0}^{\infty}t^{\lambda}G_{\bM}(1/t)dt.$$
\end{defin}

From Lemma~\ref{propiedadese} we have that the function $m$ is well defined in $\{\hbox{Re}(\lambda)\ge0\}$ and defines a continuous function in this set, and holomorphic in $\{\hbox{Re}(\lambda)>0\}$. Moreover, $m(x)$ is positive real for every $x\ge0$, so we can state the next
\begin{defin}
Let $\bM$ be a strongly regular sequence and let the function $e_{\bM}$ be constructed as in (\ref{funcione}). The sequence of positive real numbers $\mathfrak{m}=(m(p))_{p\in\N_0}$, is known as the \textit{sequence of moments} associated to $\bM$ (or to $e_{\bM}$).
\end{defin}
\begin{prop}\label{mequivm}
Let $\bM=(M_p)_{p\in\N_0}$ be a strongly regular sequence and $\mathfrak{m}=(m(p))_{p\in\N_0}$ the sequence of moments associated to $\bM$. Then $\bM$ and $\mathfrak{m}$ are equivalent.
\end{prop}
\begin{proof1}
\ We recall that $(m_p)_{p\in\N_0}$ is the sequence of quotients of $\bM$.
Firstly, we prove the existence of positive constants $C_1,C_2$ such that
\begin{equation}\label{porarriba}
m(p)\le C_1 C_2^pM_p,\qquad p\in\N_0.
\end{equation}
Let $p\in\N_0$. From Proposition~\ref{gm}$.(i)$, there exists $k_3>0$ with
$$m(p)\le\int_{0}^{\infty}t^{p}h_{\bM}(k_3/t)dt=
\int_{0}^{m_{p+1}}t^{p}h_{\bM}(k_3/t)dt+
\int_{m_{p+1}}^{\infty}t^{p}h_{\bM}(k_3/t)dt.$$
In the first integral we take into account that $h_{\bM}$ is bounded by 1, while in the second one we use the definition of $h_{\bM}$.
This yields
$$m(p)\le\int_{0}^{m_{p+1}}t^{p}dt+ \int_{m_{p+1}}^{\infty}t^{p}\frac{k_3^{p+2}}{t^{p+2}}M_{p+2}dt=
\frac{1}{p+1}m_{p+1}^{p+1}+k_3^{p+2}\frac{M_{p+2}}{m_{p+1}}.$$
We have $M_{p+2}=m_{p+1}M_{p+1}$, and we may apply the property $(\mu)$ of $\bM$ and (\ref{e107}) to obtain that
$$
m(p)\le A^2M_{p+1}+k_3^{p+2}M_{p+1}\le (A^3M_1A^p+AM_1k_3^2A^pk_3^p)M_p,
$$
as desired. This concludes the first part of the proof.\par
We will now show the existence of constants $C_3,C_4>0$ such that $m(p)\ge C_3 C_4^pM_p$ for every $p\in\N_0$. Let $p\in\N_0$. From Proposition~\ref{gm}$.(i)$, there exist $k_1,k_2>0$ such that
$$
m(p)\ge k_1\int_{0}^{\infty}t^{p}h_{\bM}(k_2/t)dt\ge k_1 \int_{0}^{k_2m_p}t^{p}h_{\bM}(k_2/t)dt.
$$
Since the map $t\mapsto h_{\bM}(k_2/t)$ decreases in $(0,\infty)$, we have that for every $t\in(0,k_2m_p]$,
$$
h_{\bM}(k_2/t)\ge h_{\bM}(1/m_p)=\frac{M_p}{m_p^p},
$$
hence
$$
m(p)\ge k_1\int_{0}^{k_2m_p}t^{p}\frac{M_p}{m_p^p}dt= k_1\frac{k_2^{p+1}m_p^{p+1}}{p+1}\frac{M_p}{m_p^p}=
k_1\frac{k_2^{p+1}}{p+1}m_pM_p.
$$
Now, $m_pM_p=M_{p+1}$ and $p+1\le 2^p$ for every $p\in\N_0$. By applying (\ref{propgrowth}) we finally conclude that
$$
m(p)\ge k_1k_2M_1(k_2/2)^pM_p.
$$
\end{proof1}

\begin{rem}
In the Gevrey case of order $\a>0$, $\bM_{\a}=(p!^{\a})_{p\in\N_0}$, we may choose
$$
e_{\bM_{\a}}(z)=\frac{1}{\a}z^{1/\a}\exp(-z^{1/\a}),\qquad z\in S_{\a}.
$$
Then we obtain that $m_{\a}(\lambda)=\Gamma(1+\a \lambda)$ for $\Re(\lambda)\ge 0$. Of course, the sequences $\bM_{\a}$ and $\mathfrak{m}_{\a}=(m_{\a}(p))_{p\in\N_0}$ are equivalent.
\end{rem}

\section{Right inverses for the asymptotic Borel map in ultraholomorphic classes in sectors}\label{sectRightInver1var}

The proof of the incoming result follows the same lines as the original one in the Gevrey case (see~\cite{to},~\cite{ca},~\cite[Thm.\ 4.1]{javier}). The only difficulty stems from the use of the kernel $e_{\bM}$, linked to a general sequence $\bM$ and, to a certain extent, unknown, whereas the exponential function linked to the Gevrey case is very well-known.

\begin{theo}\label{tpral}
Let $\bM=(M_{p})_{p\in\N_0}$ be a strongly regular sequence and let $S=S(d,\delta)$ be a sector with vertex at the origin and opening $0<\delta<\gamma(\bM)$. For every $(a_{p})_{p\in\N_0}\in\Lambda_{\bM}(\N_{0})$ there exists a function $f\in\mathcal{A}_{\bM}(S)$ such that $f$ admits $\hat{f}=\sum_{p\in\N_0}\frac{a_p}{p!}z^p$ as its uniform asymptotic expansion in $S_{\delta}$.
\end{theo}
\begin{proof1}
We may assume that $d=0$ without loss of generality, for the case $d\neq0$ only involves an adequate rotation.

Let $(a_{p})_{p\in\N_0}\in\Lambda_{\bM}(\N_0)$, and let $\mathfrak{m}=(m(p))_{p\in\N_0}$ be the sequence of moments associated to~ $\bM$.
There exist positive constants $C_1,A_1$ such that
\begin{equation}\label{equatBoundsLambdaMA}
|a_p|\le C_1D_1^{p}p!M_p,\quad p\in\N_0.
\end{equation}
From Proposition~\ref{mequivm}, the series
$$
\hat{g}=\sum_{p\in\N_0}\frac{a_p}{p!m(p)}z^p
$$
is convergent in a disc $D(0,R)$ for some $R>0$, and it defines a holomorphic function $g$ there. Let $0<R_0<R$.
We define
\begin{equation}\label{intope}
f(z):=\int_{0}^{R_0}e_{\bM}\left(\frac{u}{z}\right)g(u)\frac{du}{u},\qquad z\in S_{\delta},
\end{equation}
where the kernel $e_{\bM}$ is constructed as in (\ref{funcione}).
By virtue of Leibnitz's theorem on analyticity of parametric integrals and the definition of $e_{\bM}$, $f$ turns out to be a holomorphic function in $S_{\delta}$. Let us prove that
$f\sim_{\bM}\hat{f}$ uniformly in $S_{\delta}$.\par
Let $N\in\N$ and $z\in S_{\delta}$. We have
\begin{align*}
f(z)-\sum_{p=0}^{N-1}a_p\frac{z^p}{p!} &= f(z)-\sum_{p=0}^{N-1}\frac{a_p}{m(p)}m(p)\frac{z^p}{p!}\\
&= \int_{0}^{R_0}e_{\bM}\left(\frac{u}{z}\right)\sum_{k=0}^{\infty}\frac{a_{k}}{m(k)}\frac{u^k}{k!}\frac{du}{u} -\sum_{p=0}^{N-1}\frac{a_p}{m(p)}\int_{0}^{\infty}u^{p-1}e_{\bM}(u)du\frac{z^p}{p!}.
\end{align*}
After a change of variable $v=zu$ in the second integral, by virtue of the estimate~(\ref{propiedadese}) one may use Cauchy's residue theorem in order to check that
$$
z^p\int_{0}^{\infty}u^{p-1}e_{\bM}(u)du= \int_{0}^{\infty}v^{p}e_{\bM}\left(\frac{v}{z}\right)\frac{dv}{v},
$$
which allows us to write the preceding difference as
\begin{multline*}
\int_{0}^{R_0}e_{\bM}\left(\frac{u}{z}\right)\sum_{k=0}^{\infty}\frac{a_{k}}{m(k)}\frac{u^k}{k!}\frac{du}{u} -\sum_{p=0}^{N-1}\frac{a_p}{m(p)}\int_{0}^{\infty}u^{p}e_{\bM}\left(\frac{u}{z}\right)\frac{du}{u}\frac{1}{p!}\\
=\int_{0}^{R_0}e_{\bM}\left(\frac{u}{z}\right)\sum_{k=N}^{\infty}\frac{a_{k}}{m(k)}\frac{u^k}{k!}\frac{du}{u} -\int_{R_0}^{\infty}e_{\bM}\left(\frac{u}{z}\right)\sum_{p=0}^{N-1}\frac{a_p}{m(p)}\frac{u^{p}}{p!}\frac{du}{u}.
\end{multline*}
Then, we have
$$\left|f(z)-\sum_{p=0}^{N-1}a_p\frac{z^p}{p!}\right|\le f_{1}(z)+f_2(z),$$
where
$$f_{1}(z)=\left|\int_{0}^{R_0}e_{\bM}\left(\frac{u}{z}\right)\sum_{k=N}^{\infty}\frac{a_{k}}{m(k)}\frac{u^k}{k!}\frac{du}{u}\right|,$$
$$f_{2}(z)=\left|\int_{R_0}^{\infty}e_{\bM}\left(\frac{u}{z}\right)\sum_{p=0}^{N-1}\frac{a_p}{m(p)}\frac{u^{p}}{p!}\frac{du}{u}\right|.$$
We now give suitable estimates for $f_1(z)$ and $f_2(z)$.
From Proposition~\ref{mequivm} there exist $C_2,D_2>0$ (not depending on $z$) such that
\begin{equation}\label{e327}
\frac{a_{k}}{m(k)k!}\le \frac{C_1D_1^{k}k!M_k}{m(k)k!}\le C_2D_2^k,
\end{equation}
for all $k\in\N_0$. This yields
$$f_{1}(z)\le C_2\int_{0}^{R_0}\left|e_{\bM}\left(\frac{u}{z}\right)\right|\sum_{k=N}^{\infty}(D_2u)^{k}\frac{du}{u}.$$
Taking $R_0\le(1-\epsilon)/D_2$ for some $\epsilon>0$ if necessary, we get
$$f_1(z)\le \epsilon C_2D_2^{n}\int_{0}^{R_0}\left|e_{\bM}\left(\frac{u}{z}\right)\right|u^{N-1}du.$$
By a double application of $(i)$ in Proposition~\ref{gm} we derive
$$\left|G_{\bM}\left(\frac{z}{u}\right)\right|\le h_{\bM}\left(\frac{k_3|z|}{|u|}\right)=h_{\bM}\left(k_2\frac{k_3|z|}{k_2 u}\right)\le\frac{1}{k_1}G_{\bM}\left(\frac{k_3|z|}{k_2u}\right),$$
for some positive constants $k_1,k_2,k_3$. This yields
\begin{multline}\label{e337}
\int_{0}^{R_0}\left|e_{\bM}\left(\frac{u}{z}\right)\right|u^{N-1}du
\le\frac{1}{k_1}\int_{0}^{\infty}\left|\frac{u}{z}\right| G_{\bM}\left(\frac{k_3|z|}{k_2u}\right)u^{N-1}du\\
=\frac{1}{k_1}\int_{0}^{\infty}\frac{k_3t}{k_2} G_{\bM}\left(\frac{1}{t}\right)\left(\frac{k_3|z|t}{k_2}\right)^{N-1}\frac{k_3}{k_2}|z|dt\\
=\left(\frac{k_3}{k_2}\right)^{N+1}\frac{1}{k_{1}}|z|^{N} \int_{0}^{\infty}t^{N}G_{\bM}\left(\frac{1}{t}\right)dt=C_3D_3^{N}m(N)|z|^{N},
\end{multline}
for some $C_3,D_3>0$. The conclusion for $f_1$ is achieved from Proposition~\ref{mequivm}.
It only rests to estimate $f_{2}(z)$. We have $u^p\le R_0^pu^N/R_0^N$ for $u\ge R_0$ and $0\le p\le N-1$. So, according to~(\ref{e327}), we may write
$$\sum_{p=0}^{N-1}\frac{a_pu^p}{m(p)p!}\le\sum_{p=0}^{N-1}\frac{C_1D_1^pp!M_pu^p}{m(p)p!} \le\sum_{p=0}^{N-1}C_1D_1^pC_2D_2^pu^p\le\frac{u^N}{R_0^N}\sum_{p=0}^{N-1}C_1D_1^pC_2D_2^pR_0^p\le C_5D_5^Nu^N,$$
for some positive constants $C_5,D_5$. Then, we conclude
$$f_2(z)\le C_5D^N_5\int_{R_0}^{\infty}\left|e_{\bM}\left(\frac{u}{z}\right)\right|u^{N-1}du.$$
We come up to the end of the proof following similar estimates as in (\ref{e337}).
\end{proof1}

\begin{rem}
Given $\delta$ with $0<\delta<\gamma(\bM)$, choose $\delta_1$ such that $\delta<\delta_1<\gamma(\bM)$ and put $S_1=S(d,\delta_1)$.
For $A>0$ and for every $\ba=(a_{p})_{p\in\N_0}\in\Lambda_{\bM,A}(\N_{0})$, we have the estimates~(\ref{equatBoundsLambdaMA}) with $C_1=|\ba|_{\bM,A}$ and $D_1=A$. Since the previous result is valid in $S_1$, we obtain a function $f\in\mathcal{A}_{\bM}(S_1)$ that admits $\hat{f}=\sum_{p\in\N_0}\frac{a_p}{p!}z^p$ as its uniform asymptotic expansion in $S_{1}$. Moreover,
by taking into account in detail the way constants are modified in the course of the proof of Theorem~\ref{tpral}, one observes that there exist constants $C,D>0$, not depending on $f$, such that for every $N\in\N_0$ one has
\begin{equation}\left|f(z)-\sum_{p=0}^{N-1}\frac{a_p}{p!}z^p \right|\le (CC_1)(DD_1)^NM_{N}=
C|\ba|_{\bM,A}(DA)^NM_N,\qquad z\in S_1.\label{eqDesarAsintS1}
\end{equation}
According to Remark~\ref{remaCarlclassasympexpan}, there exists a constant $c=c(S,S_1)>0$ such that the restriction to $S$ of $f$ belongs to $\mathcal{A}_{\bM,cDA}(S)$, and moreover, from~(\ref{eqDesarAsintS1}) we get $\Vert f\Vert_{\bM,cDA,S}\le C|\ba|_{\bM,A}$. So, we have re-proved the following theorem of V. Thilliez.
\end{rem}

\begin{theo}\label{corolpral}
Under the hypotheses of Theorem~\ref{tpral}, there exists a positive constant $c\ge1$ such that for any $A>0$, the integral operator
$$ T_{\bM,A}:\Lambda_{\bM,A}(\N_{0})\longrightarrow\mathcal{A}_{\bM,cA}(S_{\delta})$$
defined in (\ref{intope}) by
$$T_{\bM,A}(\ba=(a_{p})_{p\in\N_0}):=\int_{0}^{R_0}e_{\bM}(u/z) \Big(\sum_{p=0}^{\infty}\frac{a_{p}}{m(p)}\frac{u^p}{p!}\Big)\frac{du}{u}$$
is linear and continuous and it turns out to be a right inverse for the asymptotic Borel map $\mathcal{B}$.
\end{theo}

\section{An application to the several variable setting}\label{seccSeverVaria}

As an application of the previous result, we will obtain a different construction of continuous extension operators in Carleman ultraholomorphic classes in polysectors of $\mathcal{R}^n$, obtained in~\cite{lastrasanz} by the first and the third authors as a generalization of V. Thilliez's result (see~\cite[Thm.\ 3.2.1]{thilliez}).
It is worth saying that the results in Section~\ref{sectRightInver1var} are also valid when the functions and sequences involved take their values in a complex Banach space $\mathbb{B}$. This will be crucial in the ongoing section.

Let $n\in\N$, $n\ge 2$, and fix a sequence $\bM=(M_{p})_{p\in\N_0}$ of positive real numbers. For a polysector $S$ in $\mathcal{R}^{n}$, the space $\mathcal{A}_{\bM}(S,\mathbb{B})$ consists of the holomorphic functions $f:S\to(\mathbb{B},\left\|\cdot\right\|_{\mathbb{B}})$ such that there exists $A>0$ (depending on $f$) with
\begin{equation}\label{e360}
\left\|f\right\|^{\mathbb{B}}_{\bM,A,S}:=\sup_{\bJ\in\N_{0}^{n},\bz\in S}\frac{\left\|D^{\bJ}f(\bz)\right\|_{\mathbb{B}}}{A^jj!M_j}<\infty
\end{equation}
(the notations adopted in Subsection~\ref{notation} are being applied).\par
For fixed $A>0$, $\mathcal{A}_{\bM,A}(S,\mathbb{B})$ consists of the elements in $\mathcal{A}_{\bM}(S,\mathbb{B})$ such that (\ref{e360}) holds, and the norm $\left\|\,\cdot\,\right\|^{\mathbb{B}}_{\bM,A,S}$ makes it a Banach space.
The space $\Lambda_{\bM,A}(\N_{0}^{n},\mathbb{B})$ consists of the multi-sequences $\lambda=(\lambda_{\bJ})_{\bJ\in\N_{0}^{n}}\in\N_{0}^{\mathbb{B}}$ such that
$$|\lambda|_{\bM,A,\mathbb{B}}:= \sup_{\bJ\in\N_{0}^{n}}\frac{\left\|\lambda_{\bJ}\right\|_{\mathbb{B}}}{A^{j}j!M_j}<\infty,$$
and $(\Lambda_{\bM,A}(\N_{0}^n,\mathbb{B}),|\cdot|_{\bM,A,\mathbb{B}})$ is a Banach space.\par
The elements in $\mathcal{A}_{\bM}(S,\mathbb{B})$ admit strong asymptotic expansion in $S$ as defined by H. Majima (see~\cite{Majima1,Majima2}), since this fact amounts, as shown by J. A. Hern\'andez~\cite{hernandez}, to having bounded derivatives in every subpolysector $T\ll S$.
The following facts, stated here without proof, can be found in detail in the three previous references and in~\cite{galindosanz,hernandezsanz,lastrasanz}.
The asymptotic information for such a function $f$ is given by the family
$$\mathrm{TA}(f)=\left\{f_{\balpha_{\bJ}}:\emptyset\neq \bJ\subseteq\mathcal{N},\balpha_{\bJ}\in\N_{0}^{\bJ}\right\},$$
where for every nonempty subset $\bJ$ of $\mathcal{N}$ and every $\balpha_{\bJ}\in\N_{0}^{\bJ}$, $f_{\balpha_{\bJ}}$ is defined as
$$f_{\balpha_{\bJ}}(\bz_{\bJ'})=\lim_{\bz_{\bJ}\to\textbf{0}_{\bJ}} D^{(\balpha_{\bJ},\textbf{0}_{\bJ'})}f(\bz),\quad \bz_{\bJ'}\in S_{\bJ'},$$
the limit being uniform on $S_{\bJ'}$ whenever $\bJ\neq\mathcal{N}$. This implies that $f_{\balpha_{\bJ}}\in\mathcal{A}_{\bM}(S_{\bJ'},\mathbb{B})$ (we agree that $\mathcal{A}_{\bM}(S_{\mathcal{N}'},\mathbb{B})$ is meant to be $\mathbb{B}$).
\begin{prop}[Coherence conditions]
Let $f\in\mathcal{A}_{\bM}(S,\mathbb{B})$ and
$$\mathrm{TA}(f)=\left\{f_{\balpha_{J}}:\emptyset\neq \bJ\subseteq\mathcal{N},\balpha_{\bJ}\in\N_{0}^{\bJ}\right\}.$$
Then, for every pair of nonempty disjoint subsets $\bJ$ and $\bL$ of $\mathcal{N}$, every $\balpha_{\bJ}\in\N_{0}^{\bJ}$ and $\balpha_{\bL}\in\N_{0}^{\bL}$, we have
\begin{equation}\label{e380}
\lim_{\bz_{\bL}\to\textbf{0}}D^{(\balpha_{\bL},\textbf{0}_{(\bJ\cup \bL)'})}f_{\balpha_{\bJ}}(\bz_{\bJ'})=f_{(\balpha_{\bJ},\balpha_{\bL})}(\bz_{(\bJ\cup \bL)'});
\end{equation}
the limit is uniform in $S_{(\bJ\cup \bL)'}$ whenever $\bJ\cup \bL\neq \mathcal{N}$.
\end{prop}
\begin{defin}\label{deffirst}
We say a family
$$\mathcal{F}=\left\{f_{\balpha_{\bJ}}\in\mathcal{A}_{\bM}(S_{\bJ'},\mathbb{B}):\emptyset\neq \bJ\subseteq\mathcal{N},\balpha_{\bJ}\in\N_{0}^{\bJ}\right\}$$
is \textit{coherent} if it fulfills the conditions given in (\ref{e380}).
\end{defin}
\begin{defin}
Let $f\in\mathcal{A}(S,\mathbb{B})$. The \textit{first order family} associated to $f$ is given by
$$\mathcal{B}_{1}(f):=\left\{f_{m_{\{j\}}}\in\mathcal{A}_{\bM}(S_{j'},\mathbb{B}): j\in\mathcal{N},m\in\N_{0}\right\}\subseteq \mathrm{TA}(f).$$
\end{defin}
The first order family consists of the elements in the total family that depend on $n-1$ variables. For the sake of simplicity, we will write $f_{jm}$ instead of $f_{m_{\{j\}}}$, $j\in\mathcal{N}$, $m\in\N_{0}$. As it can be seen in~\cite[Section 4]{galindosanz}, knowing $\mathcal{B}_{1}(f)$ amounts to knowing $\mathrm{TA}(f)$, and moreover, $\mathcal{B}_{1}(f)$ verifies what we call first order coherence conditions, emanating from the ones for $\mathrm{TA}(f)$.
In fact, there is a bijective correspondence between the set of coherent families (see Definition~\ref{deffirst}) and the one of coherent first order families
$$\mathcal{F}_{1}=\left\{f_{jm}\in\mathcal{A}_{\bM}(S_{j'},\mathbb{B}):j\in\mathcal{N},m\in\N_{0}\right\}.$$
\begin{defin}\label{defiFamiFirstOrder}
Let $\bM=(M_{p})_{p\in\N_0}$ be a sequence that fulfills property $(\mu)$ for a constant $A_1$, and let $A>0$. We define $\mathfrak{F}^1_{\bM,A}(S,\mathbb{B})$ as the set of coherent families of first order
$$\mathcal{G}=\left\{f_{jm}\in\mathcal{A}_{\bM,2AA_{1}}(S_{j'},\mathbb{B}):j\in\mathcal{N},m\in\N_{0}\right\}$$
such that for every $j\in\mathcal{N}$ we have
$$\mathcal{G}_{j}:=(f_{jm})_{m\in\N_0}\in \Lambda_{\bM,2AA_1}\left(\N_0,\mathcal{A}_{\bM,2AA_1}(S_{j'},\mathbb{B})\right).$$
\end{defin}
It is immediate to prove that, if we put
$$\nu_{\bM,A}(\mathcal{G}):= \sup_{j\in\mathcal{N}}\left\{ |\mathcal{G}_{j}|_{\bM,2AA_1,\mathcal{A}_{\bM,2AA_1}(S_{j'},\mathbb{B})}\right\}, \quad\mathcal{G}\in\mathfrak{F}_{\bM,A}^1(S,\mathbb{B}),$$
then $(\mathfrak{F}_{\bM,A}^1(S,\mathbb{B}),\nu_{\bM,A})$ is a Banach space.
We may consider a generalized Borel map, say $\mathcal{B}_{1}$, sending any function in $\mathcal{A}_{\bM,A}(S,\mathbb{B})$ to its corresponding first order family. Then, one has
\begin{prop}[\cite{lastrasanz},\ Proposition~3.4]\label{propBunowelldefined}
The map
$\mathcal{B}_{1}:\mathcal{A}_{\bM,A}(S,\mathbb{B})\to \mathfrak{F}_{\bM,A}^{1}(S,\mathbb{B})$
is well defined, linear and continuous.
\end{prop}
The main purpose of the current section is to obtain a continuous right inverse for the preceding operator. The procedure followed is similar to the one in \cite{javier} for Gevrey classes, and it is based on our new proof of Theorem~\ref{corolpral}, so overcoming the technical difficulties encountered in~\cite{lastrasanz}.

The first step in the proof consists of changing the problem into an equivalent one in terms of functions in one variable with values in an appropriate Banach space of functions. In order to do this, the following result is essential. Since the proof for a similar statement can be found in detail in \cite[Thm.\ 4.5]{hernandezsanz}, we omit it.

\begin{theo}\label{tcambio}
Let $n,m\in\N$, $\bM$ be a sequence of positive real numbers, $\mathbb{B}$ be a complex Banach space, $A>0$ and $S$ and $V$ be (poly)sectors in $\mathcal{R}^{n}$ and $\mathcal{R}^{m}$, respectively. Then, we have:
\begin{itemize}
\item[(i)] If $\bM$ fulfills $(\mu)$ and $A_1$ is the constant involved in this property, then the map
$$\psi_{1}:\mathcal{A}_{\bM,A}(S\times V,\mathbb{B})\longrightarrow\mathcal{A}_{\bM,2AA_{1}}\left(S,\mathcal{A}_{\bM,2AA_{1}}(V,\mathbb{B})\right)$$
sending each function $f\in \mathcal{A}_{\bM,A}(S\times V,\mathbb{B})$ to the function $f^{\star}=\psi_{1}(f)$ given by
$$\left(f^{\star}(\bz)\right)(\bw)=f(\bz,\bw),\qquad (\bz,\bw)\in S\times V,$$
is well defined, linear and continuous. Given $f\in\mathcal{A}_{\bM,A}(S\times V,\mathbb{B})$, for every $\balpha\in\N_{0}^{n}$, $\bbeta\in\N_{0}^{m}$ and $(\bz,\bw)\in S\times V$ we have
$$D^{(\balpha,\bbeta)}f(\bz,\bw)=D^{\balpha}\left(D^{\balpha}f^{\star}(\bz)\right)(\bw),$$
and so
$$\left\|f^{\star}\right\|_{\bM,2AA_1,S}^{\mathcal{A}_{\bM,2AA_{1}}(V,\mathbb{B})}\le \left\|f\right\|_{\bM,A,S\times V}^{\mathbb{B}}.$$
\item[(ii)] If $\bM$ fulfills $(\alpha_0)$, the map
$$\psi_{2}:\mathcal{A}_{\bM,A}\left(S,\mathcal{A}_{\bM,A}(V,\mathbb{B})\right) \longrightarrow \mathcal{A}_{\bM,A}(S\times V,\mathbb{B})$$
given by
$$\left(\psi_2(f)\right)(\bw,\bz)=\left(f(\bz)\right)(\bw),\qquad (\bz,\bw)\in S\times V,$$
is well defined, linear and continuous. For $f\in \mathcal{A}_{\bM,A}\left(S,\mathcal{A}_{\bM,A}(V,\mathbb{B})\right)$, every $\balpha\in\N_{0}^{n}$, $\bbeta\in\N_{0}^{m}$ and $(\bz,\bw)\in S\times V$ we have
\begin{equation}\label{igualDerivIsomor}
D^{(\balpha,\bbeta)}\left(\psi_{2}(f)\right)(\bz,\bw)=D^{\bbeta}\left(D^{\balpha}f(\bz)\right)(\bw),
\end{equation}
and consequently
$$\left\|\psi_{2}(f)\right\|_{\bM,A,S\times V}^{\mathbb{B}}\le\left\|f\right\|_{\bM,A,S}^{\mathcal{A}_{\bM,A}(V,\mathbb{B})}.$$
\end{itemize}
\end{theo}
Before stating the main result in this section, we need some information about the asymptotic behaviour of the one-variable solution provided by Theorem~\ref{corolpral} when it takes its values in a Banach space of the type $\mathcal{A}_{\bM,A}(S,\mathbb{B})$.\par

Let $n\ge 1$, $A>0$, $\bM=(M_{p})_{p\in\N_0}$ be a strongly regular sequence,
$S$ a polysector in $\mathcal{R}^n$ and $0<\delta<\gamma(\bM)$. Suppose that for every $p\in\N_{0}$ we are given a function $f_{p}\in\mathcal{A}_{\bM,A}(S,\mathbb{B})$ in such a way that $\bfe=(f_{p})_{p\in\N_0}\in\Lambda_{\bM,A}\left(\N_0,\mathcal{A}_{\bM,A}(S,\mathbb{B})\right)$. Let $R_0$ be as in the proof of Theorem~\ref{tpral}. By Theorem~\ref{corolpral}, we know that the function $H^{\star}:=T_{\bM,A}(\bfe):S_{\delta}\to\mathcal{A}_{\bM,A}(S,\mathbb{B})$, given by
$$H^{\star}(w)=\int_{0}^{R_0}e_{\bM}(u/w) \Big(\sum_{p=0}^{\infty}\frac{f_{p}}{m(p)}\frac{u^{p}}{p!}\Big)\frac{du}{u},$$
belongs to $\mathcal{A}_{\bM,c(\delta)A}(S_{\delta},\mathcal{A}_{\bM,A}(S,\mathbb{B}))$, for suitable $c(\delta)>1$, and it admits $\sum_{p\ge0}f_{p}z^p/p!$ as uniform $\bM$-asymptotic expansion in $S_{\delta}$. Hence, the function $H:S_{\delta}\times S\to\mathbb{B}$ given by $H(w,\bz)=H^{\star}(w)(\bz)$, belongs, by Theorem~\ref{tcambio}$.(ii)$, to $\mathcal{A}_{\bM,A}(S_{\delta}\times S,\mathbb{B})$ and, for every $\balpha\in\N_{0}^{n}$, we have
\begin{equation}
D^{(0,\balpha)}H(w,\bz)=D^{\balpha}(H^{\star}(w))(\bz)
=\int_{0}^{R_0}e_{\bM}(u/w) \Big(\sum_{p=0}^{\infty}\frac{D^{\balpha}f_{p}(\bz)}{m(p)}\frac{u^{p}}{p!}\Big)\frac{du}{u}.
\label{e473}
\end{equation}
The proof of the next Lemma, extremely lengthy and awkward when following the technique in~\cite{thilliez}, is now easy due to the new solution in integral form for Theorem~\ref{corolpral}.
\begin{lemma}\label{lemaDerivNula}
Let $S=\prod_{j\in\mathcal{N}}S_j$. If for every $m,p\in\N_0$ and $j\in\mathcal{N}$, we have
\begin{equation}\label{e480}
\lim_{z_{j}\to0,z_{j}\in S_{j}}D^{m\be_{j}}f_{p}(\bz)=0\hbox{ \ uniformly on }S_{j'},
\end{equation}
then, for every $m\in\N_0$ and $j\in\mathcal{N}$ one has
$$\lim_{z_{j}\to0,z_{j}\in S_{j}}D^{(0,m\be_{j})}H(w,\bz)=0\hbox{ \ uniformly on }S_{\delta}\times S_{j'}.
$$
\end{lemma}
\begin{proof1}
By (\ref{e473}) we have
\begin{equation*}
D^{(0,m\be_{j})}H(w,\bz)= \int_{0}^{R_0}e_{\bM}(u/w) \Big(\sum_{p=0}^{\infty}\frac{D^{m\be_{j}}f_{p}(\bz)}{m(p)}\frac{u^{p}}{p!}\Big)\frac{du}{u}.
\end{equation*}
Given $\varepsilon>0$, there exists $p_{0}\in\N_0$ such that, for every $p\ge p_0$, every $\bz\in S$ and every $u\in [0, R_0]$, one has
$$
\left\|\sum_{p=p_0}^{\infty}\frac{D^{m\be_{j}}f_{p}(\bz)}{m(p)}\frac{u^{p}}{p!}\right\|<\varepsilon.$$
From (\ref{e480}), there exists $M>0$ such that whenever $\bz=(z_1,...,z_n)\in S$ and $z_{j}\in S_{j}\cap D(0,M)$ we have
$$\left|D^{m\be_{j}}f_{p}(\bz)\right|\le\frac{\varepsilon m(p)p!}{p_0R_0^{p}},\qquad p=0,\ldots,p_0-1.$$
Hence, for every $\bz\in S$ with $z_{j}\in S_{j}\cap D(0,M)$ and every $w\in S_{\delta}$ we have
\begin{align*}
\left|	D^{(0,m\be_j)}H(w,\bz)\right|&\le \int_{0}^{R_0}|e_{\bM}(u/w)|\Big(\Big|\sum_{p=0}^{p_0-1}\frac{D^{m\be_{j}}f_{p}(\bz)}{m(p)} \frac{u^{p}}{p!}\Big|+\Big|\sum_{p=p_0}^{\infty}\frac{D^{m\be_{j}}f_{p}(\bz)}{m(p)} \frac{u^{p}}{p!}\Big|\Big)\frac{du}{u}\\
&\le 2\varepsilon\int_{0}^{R_0}|e_{\bM}(u/w)|\frac{du}{u}=
2\varepsilon\int_0^{\infty}\frac{1}{|w|}|G_{\bM}(w/u)|du.
\end{align*}
We will be done if the last integral is uniformly bounded in $S_{\delta}$. In the following estimates we use Proposition~\ref{gm}.$(i)$, the fact that $h_{\bM}$ is bounded above by 1, and the very definition of $h_{\bM}$:
\begin{align*}
\int_0^{\infty}\frac{1}{|w|}|G_{\bM}(w/u)|du&= \int_0^{|w|}\frac{1}{|w|}du+\int_{|w|}^{\infty}\frac{1}{|w|}h_{\bM}(k_3|w|/u)du\\
&\le 1+\int_{|w|}^{\infty}\frac{1}{|w|}k_3^2\frac{|w|^2}{u^2}M_2du=1+k_3^2M_2.
\end{align*}
\end{proof1}

We now state our extension result. The proof is partially included, for it follows similar steps as in Theorem~3.4 in \cite{javier}, or Theorem 3.6 in~\cite{lastrasanz}.
\begin{theo}
Let $\bM=(M_{p})_{p\in\N_0}$ be a strongly regular sequence and $\bdelta=(\delta_{1},...,\delta_{n})\in(0,\infty)^n$ with $\delta_{j}<\gamma(\bM)$ for $j\in\mathcal{N}$. Then, there exists a constant $c=c(\bM,\bdelta)>1$, a constant $C=C(\bM,\bdelta)>0$, and for every $A>0$, a linear operator
$$U_{\bM,A,\bdelta}:\mathfrak{F}_{\bM,A}^{1}(S_{\bdelta})\longrightarrow \mathcal{A}_{\bM,cA}(S_{\bdelta})$$
such that, for every $\mathcal{G}\in\mathfrak{F}_{\bM,A}^{1}(S_{\bdelta})$ we have
$$\mathcal{B}_{1}\left(U_{\bM,A,\bdelta}(\mathcal{G})\right)=\mathcal{G}\quad\hbox{ and }\quad \left\|U_{\bM,A,\bdelta}(\mathcal{G})\right\|_{\bM,cA,S_{\bdelta}}  \le C\nu_{\bM,A}(\mathcal{G}).$$
\end{theo}
\begin{proof1}
Suppose $\bM$ verifies $(\mu)$ for a constant $A_1>0$. Let
$$\mathcal{G}= \left\{f_{jm}\in\mathcal{A}_{\bM,2AA_{1}}(S_{\bdelta_{j'}}): j\in\mathcal{N},m\in\N_{0}\right\}\in\mathfrak{F}_{\bM,A}^{1}(S_{\bdelta}).
$$
The proof is divided into $n$ steps, in such a way that in the $k$-th step we will obtain a function whose first order family contains the first $k$ sequences $(f_{jm})_{m\in\N_0}$, with $j\le k$.
We will only detail the first two steps.\par
Since $\mathcal{G}_{1}:=\{f_{1m}\}\in \Lambda_{\bM,2AA_1}\left(\N_0,\mathcal{A}_{\bM,2AA_1}(S_{\bdelta_{1'}})\right)$,
the vector-valued version of Theorem~\ref{corolpral} provides
constants $c_{1}\ge 1$, $C_{1}>0$ and a linear continuous operator
$$
T_{\bM,2AA_{1},\delta_{1}}: \Lambda_{\bM,2AA_{1}}(\N_0,\mathcal{A}_{\bM,2AA_{1}}(S_{\bdelta_{1'}})) \longrightarrow \mathcal{A}_{\bM,c_{1}2AA_{1}}(S_{\delta_{1}}, \mathcal{A}_{\bM,2AA_{1}}(S_{\bdelta_{1'}}))
$$
such that, if we put $H_{1}^{[1]\star}:=T_{\bM,2AA_{1},\delta_{1}}(\mathcal{G}_{1})$, then
$$
H_{1}^{[1]\star}\sim_{\bM}\sum_{m=0}^{\infty}\frac{f_{1m}}{m!}z_{1}^{m}\quad\textrm{and}\quad \|H_{1}^{[1]\star}\|_{\bM,c_{1}2AA_{1},S_{\delta_{1}}}^{\mathcal{A}_{\bM,2AA_{1}}(S_{\bdelta_{1'}})} \le C_{1}|\mathcal{G}_{1}|_{\bM,2AA_{1},\mathcal{A}_{\bM,2AA_{1}}(S_{\bdelta_{1'}})}.
$$
Since
$$
\mathcal{A}_{\bM,c_{1}2AA_{1}}(S_{\delta_{1}},\mathcal{A}_{\bM,2AA_{1}}(S_{\bdelta_{1'}}))\subseteq
\mathcal{A}_{\bM,c_{1}2AA_{1}}(S_{\delta_{1}},\mathcal{A}_{\bM,c_{1}2AA_{1}}(S_{\bdelta_{1'}}))
$$
(with the correspondent inequality for the norms), by Theorem~\ref{tcambio}.$(ii)$
we know that the function $H^{[1]}:S_{\bdelta}\to \C$ given by
$$H^{[1]}(\bz):=H_1^{[1]\star}(z_{1})(\bz_{1'}),\qquad \bz=(z_1,\bz_{1'})\in S_{\bdelta},
$$
belongs to $\mathcal{A}_{\bM,c_{1}2AA_{1}}(S_{\bdelta})$ and, moreover,
$$
\|H^{[1]}\|_{\bM,c_{1}2AA_{1},S_{\bdelta}}\le \|H^{[1]\star}_{1}\|_{\bM,c_{1}2AA_{1},S_{\delta_{1}}}^{\mathcal{A}_{\bM,2AA_{1}}(S_{\bdelta_{1'}})}.
$$
Let $\mathcal{B}_{1}(H^{[1]})=\{h_{jm}^{[1]}:j\in \mathcal{N},m\in\N_{0}\}$. For every $\bz_{1'}\in S_{\bdelta_{1'}}$ we have, by virtue of (\ref{igualDerivIsomor}),
$$
h_{1m}^{[1]}(\bz_{1'})=\lim_{z_{1}\to 0,z_{1}\in S_{\delta_{1}}}D^{m\be_{1}}H^{[1]}(\bz)= \lim_{z_{1}\to 0,z_{1}\in S_{\delta_{1}}}(H^{[1]\star}_{1})^{(m)}(z_{1})(\bz_{1'})=f_{1m}(\bz_{1'}).
$$
This concludes the first step of the proof. Let $H_{2}^{[1]\star}$ be the function given by
$$
H_{2}^{[1]\star}(z_{2})(\bz_{2'}):=H^{[1]}(z_{2},\bz_{2'}),\qquad z_{2}\in S_{\delta_{2}},\ \bz_{2'}\in S_{\bdelta_{2'}}.
$$
From Theorem~\ref{tcambio}.$(i)$, we have
$$
H_{2}^{[1]\star}\in\mathcal{A}_{\bM,c_{1}(2A_{1})^{2}A}(S_{\delta_{2}},\mathcal{A}_{\bM,c_1(2  A_{1})^{2}A}(S_{\bdelta_{2'}})).
$$
We put
$$
H_{2}^{[1]\star}\sim_{\bM}\sum_{m=0}^{\infty}\frac{h_{2m}^{[1]}}{m!}z_{2}^{m}
$$
and, for the sake of brevity, $\mathbb{B}_{2}:=\mathcal{A}_{\bM,c_{1}(2A_{1})^{2}A}(S_{\bdelta_{2'}})$.
As $H^{[1]}\in\mathcal{A}_{\bM,c_{1}2A_{1}A}(S_{\bdelta})$, Proposition~\ref{propBunowelldefined} tells us that
$$
(h_{2m}^{[1]})_{m\in\N_{0}}\in\Lambda_{\bM,c_{1}(2A_{1})^{2}A}(\N_0,\mathbb{B}_{2}),
$$
and Definition~\ref{defiFamiFirstOrder} implies
$$
\mathcal{G}_{2}:=(f_{2m})_{m\in\N_{0}}\in \Lambda_{\bM,2A_{1}A}(\N_0,\mathcal{A}_{\bM,2A_{1}A}(S_{\bdelta_{2'}})).
$$
So, $(f_{2m}-h_{2m}^{[1]})_{m\in\N_{0}}\in\Lambda_{\bM,c_{1}(2A_{1})^{2}A}(\N_0,\mathbb{B}_{2})$.
By Theorem~\ref{corolpral}, we have $c_{2}\ge 1$, $C_{2}>0$ and a linear continuous operator
$$
T_{\bM,c_{1}(2A_{1})^{2}A,\delta_{2}}:\Lambda_{\bM,c_{1}(2A_{1})^{2}A}(\N_0,\mathbb{B}_2) \longrightarrow \mathcal{A}_{\bM,c_{2}c_{1}(2A_{1})^{2}A}(S_{\delta_{2}},\mathbb{B}_2)
$$
such that, if we define
$$
H_{2}^{[2]\star}:=T_{\bM,c_{1}(2A_{1})^{2}A,\delta_{2}}\big((f_{2m}-h_{2m}^{[1]})_{m\in\N_{0}}\big),
$$
then
\begin{equation}\label{desaasinh22star}
H_{2}^{[2]\star}\sim_{\bM}\sum_{m=0}^{\infty}\frac{f_{2m}-h^{[1]}_{2m}}{m!}z_{2}^{m}
\end{equation}
and
$$
\|H_{2}^{[2]\star}\|_{\bM,c_{2}c_{1}(2A_{1})^{2}A,S_{\delta_{2}}}^{\mathbb{B}_2}\le C_{2}|(f_{2m}-h_{2m}^{[1]})_{m\in\N_{0}}|_{\bM,c_{1}(2A_{1})^{2}A,\mathbb{B}_{2}}.
$$
Since
$$
\mathcal{A}_{\bM,c_{2}c_{1}(2A_{1})^{2}A}(S_{\delta_{2}},\mathbb{B}_{2})\subseteq \mathcal{A}_{\bM,c_{2}c_{1}(2A_{1})^{2}A}(S_{\delta_{2}}, \mathcal{A}_{\bM,c_{2}c_{1}(2A_{1})^{2}A}(S_{\bdelta_{2'}})),
$$
$H_{2}^{[2]\star}$ belongs to the second of these spaces, and Theorem~\ref{tcambio}.$(ii)$ ensures that the function $H^{[2]}:S_{\bdelta}\to \C$ given by
$$H^{[2]}(\bz):=H^{[2]\star}_2(z_{2})(\bz_{2'}),\qquad \bz=(z_2,\bz_{2'})\in S_{\bdelta},
$$
belongs to $\mathcal{A}_{\bM,c_{2}c_{1}(2A_{1})^{2}A}(S_{\bdelta})$ and
$$
\|H^{[2]}\|_{\bM,c_{2}c_{1}(2A_{1})^{2}A,S_{\bdelta}}\le \|H^{[2]\star}_{2}\|_{\bM,c_2c_{1}(2A_{1})^{2}A,S_{\delta_{2}}}^{\mathbb{B}_2}.
$$
We write $\mathcal{B}_{1}(H^{[2]})=\{h_{jm}^{[2]}: j\in\mathcal{N}, m\in\N_{0}\}$. For $j=1$, due to the coherence conditions for the families $\mathcal{G}$ and $\mathcal{B}_{1}(H^{[1]})$ we have for all $m,k\in\N_{0}$,
$$
\lim_{z_{1}\to0,z_{1}\in S_{\delta_{1}}}D^{me_{1}}(f_{2k}-h^{[1]}_{2k})(\bz_{2'})=\lim_{z_{2}\to0,z_{2}\in S_{\delta_{2}}}D^{ke_{2}}(f_{1m}-h^{[1]}_{1m})(\bz_{1'})=0,$$
uniformly in $S_{\bdelta_{\{1,2\}'}}$. So, we can apply Lemma~\ref{lemaDerivNula} to guarantee that for every $m\in\N_{0}$,
we have
$$\lim_{z_{1}\to0,z_{1}\in S_{\delta_{1}}}\big(H_2^{[2]\star}\big)^{(m)}(z_1)(\bz_{1'})=0\quad\hbox{ uniformly in } S_{\bdelta_{1'}},$$
and consequently, by (\ref{igualDerivIsomor}) we deduce that for every $\bz_{1'}\in S_{\bdelta_{1'}}$,
$$
h_{1m}^{[2]}(\bz_{1'})=\lim_{z_{1}\to0,z_{1}\in S_{\delta_{1}}}D^{me_{1}}H^{[2]}(z_1,\bz_{1'})=\lim_{z_{1}\to0,z_{1}\in S_{\delta_{1}}}\big(H_2^{[2]\star}\big)^{(m)}(z_1)(\bz_{1'})=0.
$$
On the other hand, taking~(\ref{desaasinh22star}) into account, for every $\bz_{2'}\in S_{\bdelta_{2'}}$ we have
\begin{align*}
h_{2m}^{[2]}(\bz_{2'})&=\lim_{z_{2}\to0,z_{2}\in S_{\delta_{2}}}D^{me_{2}}H^{[2]}(z_2,\bz_{2'})\\
&=
\lim_{z_{2}\to0,z_{2}\in S_{\delta_{2}}}(H_{2}^{[2]\star})^{(m)}(z_{2})(\bz_{2'})=(f_{2m}-h_{2m}^{[1]})(\bz_{2'}).
\end{align*}
In conclusion, the function $F^{[2]}:=H^{[1]}+H^{[2]}$ belongs to $\mathcal{A}_{\bM,c_{2}c_{1}(2A_{1})^{2}A}(S_{\bdelta})$ and,
if we put $\mathcal{B}_{1}(F^{[2]})=\{f_{jm}^{[2]}:j\in \mathcal{N},m\in\N_{0}\}$, for every $m\in\N_{0}$ we have $f_{1m}^{[2]}=f_{1m}$, $f_{2m}^{[2]}=f_{2m}$,
and the second step is completed. We are done if $n=2$, otherwise we may repeat the previous argument until the family $\mathcal{G}$ is completely interpolated.
\end{proof1}

\section{On $\bM-$summability}

In this last section we provide some keys leading to a suitable definition of summability in general ultraholomorphic classes. First, we need to introduce some formal and analytic transforms.

\subsection{Formal and analytic $\bM-$Laplace operators}

The next definition resembles that of functions of exponential growth, playing a fundamental role when dealing with Laplace and Borel transforms in $k-$summability for Gevrey classes.
For convenience, we will say a holomorphic function $f$ in a sector $S$ is {\it continuous at the origin} if $\lim_{\bz\to 0,\ \bz\in T}f(\bz)$ exists for every $T\ll S$.

\begin{defin}
Let $\bM=(M_{p})_{p\in\N_0}$ be a strongly regular sequence, and consider a sector $S$ in $\mathcal{R}$.
The set $\mathcal{A}^{(\bM)}(S)$ consists of the holomorphic functions $f$ in $S$, continuous at 0 and such that for every unbounded proper subsector $T$ of $S$ there exist $r,k_4,k_5>0$ such that for every $z\in T$ with $|z|\ge r$ one has
\begin{equation}\label{growthhM}
|f(z)|\le\frac{k_4}{h_{\bM}(k_5/|z|)}.
\end{equation}
\end{defin}
\begin{rem}
Since continuity at 0 has been asked for, $f\in\mathcal{A}^{(\bM)}(S)$ implies that for every $T\prec S$ there exist $k_4,k_5>0$ such that for every $z\in T$ one has (\ref{growthhM}).
\end{rem}
\noindent We are ready for the introduction of the $\bM-$Laplace transform.
\begin{defin}
Let $S=S(d,\alpha)$, $f\in\mathcal{A}^{(\bM)}(S)$, $\tau\in\R$ with $|\tau-d|<\alpha\frac{\pi}{2}$ and $0<\delta<\gamma(\bM)$. Consider the function $e_{\bM}$ defined in~(\ref{funcione}). We define the \textit{$\bM-$Laplace transform of $f$ in direction $\tau$} as
$$\left(\mathcal{L}_{\bM}^{\tau}f\right)(z):= \int_{0}^{\infty(\tau)}e_{\bM}\left(\frac{u}{z}\right)f(u)\frac{du}{u},$$
for every $z\in\mathcal{R}$ with $|z|$ small enough and with $|\hbox{arg}(z)-\tau|<\delta\frac{\pi}{2}$, and where the integral is taken along the path parameterized by $t\in(0,\infty)\mapsto te^{i\tau}$.
\end{defin}

\begin{prop}\label{hollap}
Under the hypotheses of the preceding definition, $\mathcal{L}_{\bM}^{\tau}f(z)$ is well-defined and it turns out to be a holomorphic function. Moreover,
$\{\mathcal{L}_{\bM}^{\tau}f\}_{\tau/|\tau-d|<\alpha\pi/2}$ defines a holomorphic function $\mathcal{L}_{\bM}f$ in a
sectorial region $G(d,\alpha+\delta)$.
\end{prop}
\begin{proof1}
For every $u,z\in\mathcal{R}$ with $\mathrm{arg}(u)=\tau$ and $|\mathrm{arg}(z)-\tau|<\delta\frac{\pi}{2}$ we have that $u/z\in S_{\delta}$ and
$$\left|\frac{1}{u}e_{\bM}\left(\frac{u}{z}\right)f(u)\right|\le \frac{1}{|u|}\left|\frac{u}{z}\right|\left|G_{\bM}\left(\frac{z}{u}\right)\right||f(u)|\le \frac{1}{|z|}h_{\bM}\left(k_3\left|\frac{z}{u}\right|\right)\frac{k_4}{h_{\bM}(k_5/|u|)},
$$
for some positive constants $k_3,k_4,k_5$. From (\ref{e120}), the previous expression can be upper bounded by
$$\frac{k_4}{|z|}\frac{h_{\bM}^2(\rho(2)k_3|z|/|u|)}{h_{\bM}(k_5/|u|)}.$$
If we assume $|z|\le L:=k_5/(\rho(2)k_3)$, from the monotonicity of $h_{\bM}$ we derive
$$\left|\frac{1}{u}e_{\bM}\left(\frac{u}{z}\right)f(u)\right|\le \frac{k_4}{|z|}h_{\bM}\left(\frac{k_5}{|u|}\right).$$
The right part of the last inequality is an integrable function of $|u|$ in $(0,\infty)$, and Leibnitz's rule for parametric integrals allows us to conclude the first part of the proof.
Let $\sigma\in\R$ with $|\sigma-d|<\alpha\frac{\pi}{2}$. The map $\mathcal{L}^{\sigma}_{\bM}f$ is a holomorphic function in
$$\{z\in\mathcal{R}:|\hbox{arg}(z)-\sigma|<\delta\frac{\pi}{2},\ |z| \mathrm{\ small}\}.
$$
Since we know that
\begin{equation*}
\lim_{|u|\to\infty}|u|h_{\bM}\left(\frac{k_5}{|u|}\right)=0,
\end{equation*}
by Cauchy's residue theorem we easily deduce that $\mathcal{L}_{\bM}^{\tau}f(z)\equiv\mathcal{L}_{\bM}^{\sigma}f(z)$ whenever both maps are defined. Thus we can extend $\mathcal{L}_{\bM}^{\tau}f$ to a function, $\mathcal{L}_{\bM}f$, holomorphic in a sectorial region $G(d,\alpha+\delta)$.
\end{proof1}
Let $\bM=(M_{p})_{p\in\N_0}$ be a strongly regular sequence and $S=S(0,\alpha)$. It is clear that for every $\lambda\in\C$ with $\hbox{Re}(\lambda)\ge0$, the function $f_{\lambda}(z)=z^{\lambda}$ belongs to the space $\mathcal{A}^{(\bM)}(S)$. From Proposition~\ref{hollap}, one can define $\mathcal{L}_{\bM}f_{\lambda}(z)$ for every $z$ in an appropriate sectorial region $G$. Moreover, for $z\in G$ and an adequate choice of $\tau\in\R$ one has
$$\mathcal{L}_{\bM}f_{\lambda}(z)= \int_{0}^{\infty(\tau)}e_{\bM}\left(\frac{u}{z}\right)u^{\lambda-1}du.$$
In particular, for $z\in\mathcal{R}$ with $\hbox{arg}(z)=\tau$, the change of variable $u/z=t$ turns the preceding integral into
$$\int_{0}^{\infty}e_{\bM}(t)z^{\lambda-1}t^{\lambda-1}zdt=m(\lambda)z^{\lambda}.$$
Therefore, it is adequate to make the following definitions.

\begin{defin}
Given a strongly regular sequence $\bM$, the formal $\bM-$Laplace transform $\hat{\mathcal{L}}_{\bM}:\C[[z]]\to\C[[z]]$ is given by
$$\hat{\mathcal{L}}_{\bM}\left(\sum_{p=0}^{\infty}a_{p}z^{p}\right):=\sum_{p=0}^{\infty}m(p)a_{p}z^{p},\qquad \sum_{p=0}^{\infty}a_{p}z^{p}\in\C[[z]].$$
Accordingly, we define the formal $\bM-$Borel transform $\hat{\mathcal{B}}:\C[[z]]\to\C[[z]]$ by
$$\hat{\mathcal{B}}_{\bM}\left(\sum_{p=0}^{\infty}a_{p}z^{p}\right):=\sum_{p=0}^{\infty}\frac{a_{p}}{m(p)}z^{p},\qquad \sum_{p=0}^{\infty}a_{p}z^{p}\in\C[[z]].$$
\end{defin}
The operators $\hat{\mathcal{B}}_{\bM}$ and $\hat{\mathcal{L}}_{\bM}$ are inverse to each other.

\subsection{Results on quasi-analyticity in ultraholomorphic classes}

We restrict our attention to the one-variable case, although the next results are available also for functions of several variables (see~\cite{lastrasanz1}). First, quasi-analytic Carleman classes are defined.

\begin{defin}
Let $S$ be a sector in $\mathcal{R}$.
We say that $\mathcal{A}_{\bM}(S)$ is \textit{quasi-analytic} if the conditions:
\begin{description}
\item[(i)] $f\in\mathcal{A}_{\bM}(S)$, and
\item[(ii)] $\mathcal{B}(f)$ is the null sequence (or $f$ admits the null series as asymptotic expansion in $S$),
\end{description}
together imply that $f$ is the null function in $S$.
\end{defin}

Characterizations of quasi-analyticity are available for general sequences $\bM$ in \cite{lastrasanz1}, but we will focus on the case of strongly regular sequences, in which the following version of Watson's lemma may be obtained.

\begin{theo}[\cite{lastrasanz1}, Thm.\ 4.10]\label{teorWatson}
Let $\bM$ be strongly regular and let us suppose that
\begin{equation}
\sum_{n=0}^{\infty}\Big(\frac{M_{n}}{M_{n+1}}\Big)^{1/\ga(\bM)} =\infty.\label{conditionWatson}
\end{equation}
Let $\gamma\in(0,\infty)$. The following statements are equivalent:
\begin{description}
\item[(i)] $\gamma\ge\gamma(\bM)$.
\item[(ii)] The class $\mathcal{A}_{\bM}(S_{\gamma})$ is quasi-analytic.
\end{description}
\end{theo}

\begin{rem}
The proof of the implication (ii)$\Rightarrow$(i) does not need to assume (\ref{conditionWatson}). However, it is an open problem to decide whether the condition $\gamma\ge\gamma(\bM)$ implies $\mathcal{A}_{\bM}(S_{\gamma})$ is quasi-analytic without the additional assumption~(\ref{conditionWatson}), which is indeed satisfied by Gevrey sequences.
\end{rem}

\subsection{A concept of $\bM-$summability in a direction}

\noindent We finally put forward a definition of summability adapted to these general situation.

\begin{defin}
Let $\bM$ be strongly regular sequence verifying condition~(\ref{conditionWatson}). A formal power series $\hat{f}=\sum_{n\ge 0}\displaystyle\frac{f_n}{n!}z^n$ is \textit{$\bM-$summable in direction $d\in\R$} if the following conditions hold:
 \begin{description}
 \item[(i)] $(f_n)_{n\in\N_0}\in\Lambda_{\bM}$, so that $g:=
\hat{\mathcal{B}}_{\bM}\hat{f}$ converges in a disc, and
\item[(ii)] $g$ admits analytic continuation in a sector $S=S(d,\varepsilon)$ for suitable $\varepsilon>0$, and $g\in\mathcal{A}^{(\bM)}(S)$.
\end{description}
\end{defin}

\begin{prop}
Let $\hat{f}=\sum_{n\ge 0}\displaystyle\frac{f_n}{n!}z^n$ be $\bM-$summable in direction $d\in\R$.
Then, there exists a sectorial region
$G(d,\beta)$, with $\beta>\gamma(\bM)$, and a function $f\in\mathcal{A}_{\bM}(G(d,\beta))$ such that $f\sim_{\bM}\hat{f}$ in $G(d,\beta)$.
\end{prop}
\begin{proof1}
For $g$ as in the previous definition, choose $\delta>0$ such that $\delta<\gamma(\bM)<\delta+\varepsilon$, and consider the functions $G_{\bM}$ and $e_{\bM}$ defined in $S_{\delta}$. We will see that
$f:=\mathcal{L}_{\bM}g$ solves the problem. Indeed, by Proposition~\ref{hollap} we know that $\mathcal{L}_{\bM}g$ is a holomorphic function
in a sectorial region $G(d,\delta+\varepsilon)$, so that the first part of the claim is proved. In what follows we study the asymptotic behaviour of $f$. Suppose $g$ converges in the disc $D(0,R)$, and take
$0<R_0<R$. For $\tau\in\R$ with $|\tau-d|<\varepsilon\frac{\pi}{2}$ and $z\in\mathcal{R}$ with $|z|$ small enough and $|\hbox{arg}(z)-\tau|<\delta\frac{\pi}{2}$, we have
$$
f(z)=\int_{0}^{\infty(\tau)}e_{\bM}\left(\frac{u}{z}\right)g(u)\frac{du}{u}=
\int_{0}^{R_0e^{i\tau}}e_{\bM}\left(\frac{u}{z}\right)g(u)\frac{du}{u} +
\int_{R_0e^{i\tau}}^{\infty(\tau)} e_{\bM}\left(\frac{u}{z}\right)g(u)\frac{du}{u}.$$
Repeating the arguments in the proof of Theorem~\ref{tpral}, we may similarly get that
$$
\int_{0}^{R_0(\tau)}e_{\bM}\left(\frac{u}{z}\right)g(u)\frac{du}{u} \sim_{\bM}\hat{f}
$$
uniformly for $z$ as specified. On the other hand, since $g\in\mathcal{A}^{(\bM)}(S)$, there exist $k_4,k_5>0$ such that
$$
|g(u)|\le\frac{k_4}{h_{\bM}(k_5/|u|)}
$$
for every $u$ with $\mathrm{arg}(u)=\tau$. Taking into account the definition of $e_{\bM}$ and Proposition~\ref{gm}.$(i)$, there exists $k_3>0$ such that, for $z$ as before,
\begin{equation}\label{boundG2_1}
\left|\int_{R_0e^{i\tau}}^{\infty(\tau)} e_{\bM}\left(\frac{u}{z}\right)g(u)\frac{du}{u}\right|\le
\int_{R_0}^{\infty}\frac{1}{|z|}h_{\bM}\big(k_3|z|/u\big)
\frac{k_4}{h_{\bM}(k_5/u)}\,du.
\end{equation}
Now, we apply (\ref{e120}) for $s=2$, and observe that, since $h_{\bM}$ is non-decreasing, whenever $|z|<k_5/(k_3\rho(2))$ one has
\begin{equation}\label{boundG2_2}
\frac{h_{\bM}\big(k_3|z|/u\big)}{h_{\bM}(k_5/u)}\le
\frac{h^2_{\bM}\big(k_3\rho(2)|z|/u\big)}{h_{\bM}(k_5/u)}\le
h_{\bM}\big(k_3\rho(2)|z|/u\big).
\end{equation}
By Proposition~\ref{gm}.$(i)$, there exist $k_1,k_2>0$ such that
\begin{equation}\label{boundG2_3}
h_{\bM}\big(k_3\rho(2)|z|/u\big)\le
\frac{1}{k_1}G_{\bM}\Big(\frac{k_3\rho(2)|z|}{k_2u}\Big).
\end{equation}
Also, for $u>R_0$ and every $n\in\N_0$ we have $1\le (u/R_0)^n$.
So, gathering (\ref{boundG2_2}) and (\ref{boundG2_3}), the right hand side in (\ref{boundG2_1}) may be bounded above by
\begin{multline}\label{boundG2_4}
\frac{k_4}{R_0^n|z|k_1}\int_{R_0}^{\infty}u^n G_{\bM}\Big(\frac{k_3\rho(2)|z|}{k_2u}\Big)\,du
=\frac{k_4}{R_0^n|z|k_1}\Big(\frac{k_3\rho(2)|z|}{k_2}\Big)^{n+1} \int_{0}^{\infty}t^nG_{\bM}(1/t)\,dt\\
=\frac{k_4k_3\rho(2)}{k_1k_2}\Big(\frac{k_3\rho(2)}{k_2R_0}\Big)^{n}m(n)|z|^n.
\end{multline}
Since $\mathfrak{m}$ and $\bM$ are equivalent, from (\ref{boundG2_4}) we deduce that
$$
\int_{R_0e^{i\tau}}^{\infty(\tau)} e_{\bM}\left(\frac{u}{z}\right)g(u)\frac{du}{u}\sim_{\bM}\hat{0}
$$
uniformly in $A_\tau=\{z\in\mathcal{R}: |z|<k_5/(k_3\rho(2)),\ |\hbox{arg}(z)-\tau|<\delta\frac{\pi}{2}\}$ (where $\hat{0}$ is the null formal power series).
Since any $T\ll G(d,\delta+\varepsilon)$ of small radius may be covered by finitely many 
sets~$A_\tau$, the conclusion that $f\sim_{\bM}\hat{f}$ in $G(d,\delta+\varepsilon)$ is reached.

\end{proof1}
\begin{rem}
In the situation of the previous result, by Theorem~\ref{teorWatson} we deduce that $f$ is unique with the property  that $f\sim_{\bM}\hat{f}$ in $G(d,\beta)$, and it is called the \textit{$\bM-$sum of $\hat{f}$ in direction~$d$}.\par
The properties of this concept are currently under study, as well as its application to the study of solutions of different types of algebraic and differential equations in the complex domain with coefficients in general ultraholomorphic classes.
\end{rem}

\noindent\textbf{Acknowledgements}: The first and third authors are partially supported by the Spanish Ministry of Science and Innovation under project MTM2009-12561.
The second author is partially supported by the French ANR-10-JCJC 0105 project.

\vskip.5cm
\noindent Authors' Affiliation:\par\vskip.5cm

Alberto Lastra\par
Departamento de Matem\'aticas\par
Edificio de Ciencias. Campus universitario\par
Universidad de Alcal\'a\par
28871 Alcal\'a de Henares, Madrid, Spain\par
E-mail: alberto.lastra@uah.es\par
\vskip1cm
St\'ephane Malek\par
UFR de Math\'ematiques Pures et Appliqu\'ees\par
Cit\'e Scientifique - B\^at. M2\par
59655 Villeneuve d'Ascq Cedex, France\par
E-mail: malek@math.univ-lille1.fr\par
\vskip1cm
Javier Sanz\par
Departamento de An\'alisis Matem\'atico\par
Instituto de Investigaci\'on en Matem\'aticas de la Universidad de Valladolid, IMUVA\par
Facultad de Ciencias\par
Universidad de Valladolid\par
47005 Valladolid, Spain\par
E-mail: jsanzg@am.uva.es

\end{document}